\documentclass[conference]{IEEEtran}
\ifCLASSINFOpdf
  \usepackage[pdftex]{graphicx}
\else
  \usepackage[dvips]{graphicx}
\fi
\usepackage[cmex10]{amsmath}
\usepackage{algorithmic}
\usepackage{array}
\ifCLASSOPTIONcompsoc
  \usepackage[caption=false,font=normalsize,labelfont=sf,textfont=sf]{subfig}
\else
  \usepackage[caption=false,font=footnotesize]{subfig}
\fi
\usepackage{dblfloatfix}
\usepackage{booktabs}
\usepackage{url}
\usepackage{float}
\usepackage{multirow}
\usepackage{enumerate}
\usepackage{amsthm}

\usepackage{svg}
\usepackage{adjustbox}
\usepackage{booktabs,threeparttable}
\usepackage{array} % <-- new
\newcolumntype{C}[1]{>{\centering\arraybackslash}p{#1}} % <-- new

\newtheoremstyle{mystyle}%                % Name
  {}%                                     % Space above
  {}%                                     % Space below
  {\itshape}%                                     % Body font
  {}%                                     % Indent amount
  {\bfseries}%                            % Theorem head font
  {.}%                                    % Punctuation after theorem head
  { }%                                    % Space after theorem head, ' ', or \newline
  {}%                                     % Theorem head spec (can be left empty, meaning `normal')
\theoremstyle{mystyle}

\newenvironment{talign*}
 {\let\displaystyle\textstyle\csname align*\endcsname}
 {\endalign}
\usepackage{arydshln}
\usepackage[normalem]{ulem}
\usepackage{hyperref}
\usepackage{amsmath}
\usepackage{xurl}
\usepackage{cite}

 \IEEEoverridecommandlockouts
\begin{document}

%
% paper title
% can use linebreaks \\ within to get better formatting as desired
% Do not put math or special symbols in the title.
%\title{Explainable AI (XAI) for Detection of Adversarial Attacks on Neural Networks\\
\title{GraphDAC: A Graph-Analytic Approach to Dynamic Airspace Configuration
}

\author{
    Ke Feng$^{1a}$, Dahai Liu$^{2b}$, Yongxin Liu$^{2b}$, Hong Liu$^{2b}$, Houbing Song$^{1a}$\\
    % Wenkai Tan$^{~1a}$, Justus Renkhoff$^{~2b}$, Alvaro Velasquez$^{3c}$, Jian Wang$^{4e}$\\ Shuteng Niu$^{5f}$, Yongxin Liu $^{1d}$, Houbing Song$^{2b}$\\
    $^{1}$University of Maryland, Baltimore County, MD 21250 USA\\
    $^{2}$Embry-Riddle Aeronautical University, FL 32114 USA\\
    
    %$^{3}$University of Colorado Boulder, CO 80309 USA,
    % $^{4}$University of Tennessee at Martin, TN 38237 USA\\
    % $^{5}$Bowling Green State University, OH 43403 USA,

    $^{a}$\{kfeng1,songh\}@umbc.edu,
    % $^{c}$alvaro.velasquez@colorado.edu,
    $^{b}$\{liu89b, LIUY11, liuho\}@erau.edu

    % $^{e}$jwang186@utm.edu, $^{f}$sniu@bgsu.edu
    % $^{g}$zwang007@odu.edu, $^{h}$lusili@cs.odu.edu\\
    % \thanks{}
}

% The paper headers
\markboth{IEEE Internet of Things Journal,~Vol.~11, No.~4, May~2021}%
{Shell \MakeLowercase{\textit{et al.}}: Bare Demo of IEEEtran.cls for Journals}
\IEEEtitleabstractindextext{%
\begin{abstract}

The current National Airspace System (NAS) is reaching capacity due to increased air traffic, and is based on outdated pre-tactical planning. This study proposes a more dynamic airspace configuration (DAC) approach that could increase throughput and accommodate fluctuating traffic, ideal for emergencies. The proposed approach constructs the airspace as a constraints-embedded graph,  compresses its dimensions, and applies a spectral clustering-enabled adaptive algorithm to generate collaborative airport groups and evenly distribute workloads among them. Under various traffic conditions, our experiments demonstrate a 50\% reduction in workload imbalances. This research could ultimately form the basis for a recommendation system for optimized airspace configuration. Code available at \url{https://github.com/KeFenge2022/GraphDAC.git}.

\end{abstract}

}

\IEEEoverridecommandlockouts
\maketitle
\IEEEdisplaynontitleabstractindextext
\IEEEpeerreviewmaketitle

\section{Introduction}
% Step 1 - topic
% what is the current airspace condition? 
The National Airspace System (NAS) comprises a complex interplay of airports and facilities, ensuring safe and smooth air travel. Air traffic control (ATC), primarily handled by human controllers, issues directives to pilots to avoid collisions and other dangers \cite{gianazza2010}. In ATC, controllers' workload is crucial to safety \cite{lee2008examining} and is mitigated by manually dividing and merging airspace into groups with dedicated personnel. This practice is termed Airspace Configuration (AC). 

% limitation of AC
Airspace is typically pre-configured according to historical plans, with minor alterations by Air Traffic Control (ATC) managers \cite{zelinski2011comparing,lee2008examining,gianazza2010}. However, with growing use of advanced technologies like Unmanned Aircraft Systems (UAS), this reliance on human expertise may prove insufficient for the increasing complexity of airspace dynamics \cite{wang2021counter}. Furthermore, conventional ATC struggles to react swiftly to emergencies, causing traffic congestion and delays \cite{yang2022multiagent, yang2021machine}. As a result, Dynamic Airspace Configuration (DAC), a real-time, data-driven approach, has gained attention. Unlike static historical models, DAC adjusts to traffic demand while accommodating constraints such as weather, fleet diversity, congestion, and sector complexity \cite{kopardekar2007initial, kopardekar2008airspace}.

% Step 2 - problem 
% what is DAC? % limitation of current DAC
Several methods have been proposed for Dynamic Airspace Configuration (DAC), but their real-world effectiveness is debatable. For example, Dynamic Airspace Sectorization \cite{Sergeevagenetic} overlooks controllers' coordination workload and reconfiguration cost, creating entirely new configurations for each day's segment. SectorFlow \cite{brinton2008airspace} groups flight trajectories to minimize airspace complexity, assigning airspace to each cluster. Its improved version refines initial partition using gradient search and keeps flow intersections off sector boundaries \cite{brinton2009airspace}. CellGeoSect \cite{sabhnani2010flow}, a cell clustering method, visualizes the airspace as hexagonal cells, maximizes flow connectivity, and balances flight counts between clusters. It then modifies the design to avoid significant flow's geometric constraints. However, these methods can decrease efficiency as controllers may be unfamiliar with newly assigned airports and consequently bring safety concerns.

This study introduces a practical, adaptive algorithm for airspace configuration, which modifies existing configurations with minimal changes, rather than designing entirely new ones. This method entails a three-stage graph-based clustering method. Firstly, we built simulated airspace from open-source data. We then convert the airspace into a relation graph that embeds the operational constraints. This is done by only connecting the geographically adjacent airports and setting their edge weights negatively related to their gross workloads. Secondly, considering that airports within the relation graph are sparsely connected, we then increase the computation efficiency by mapping each airport in the relation graph into a low-dimension space, in which Singular value decomposition (SVD)\cite{von2007tutorial,galluccio2012graph} and Autoencoder\cite{wang2016auto} are compared. Finally, we perform a spectral clustering-enabled adaptive algorithm on the low-dimensional space to get the new configuration with the traffic-center pattern surrounded by the non-busy airports. The contributions of our work are as follows:
\begin{itemize}
    \item We propose a graph clustering-enabled algorithm for DAC: we allow configuration change around adjacent airports to minimize collaboration costs while balancing the ATC controllers' workload.
    \item We propose a three-stage graph clustering method: we construct a graph from airspace with embedded spatial-temporal constraints, then reduce the graph dimension and perform an adaptive clustering to get the final sector configuration, in which each busy sector is regarded as the center and is surrounded by non-busy airports. Our experiments show that our method can reduce sector workload unbalanced level by over 50\% in different traffic conditions.
    \item We investigated the efficacy of reducing the dimensions of a graph using linear and non-linear methods: SVD (Singular Value Decomposition) and Autoencoder. Our findings demonstrate that SVD effectively reduces collaboration workload when transitioning to new configurations. On the other hand, the autoencoder excels at minimizing workload imbalances for new configurations.
\end{itemize}

The remainder of this paper is organized as follows: A literature review of related work is presented in Section~\ref{sectRW}. We present the methodology in Section~\ref{sectMM}. Evaluation and discussion are presented in Section~\ref{sectEED} and conclusions in Section~\ref{sectCC}.

\section{Related Work}
\label{sectRW}

% \subsection{Dynamic Airspace Configuration}

The leading solution for Dynamic Airspace Configuration (DAC) currently involves tactical, dynamic adjustment of airspace to minimize demand and capacity imbalances \cite{zelinski2011comparing}. This section will summarize commonly employed solutions.

% First, dynamically creating entirely new controlled airspace. Controlled airspace sector boundaries are newly created every time without relying on pre-existing structures. This approach is also referred to as dynamic sectorization. For example, \textit{SectorFlow}\cite{brinton2008airspace} clusters flight trajectory and tries to minimize airspace complexity simultaneously. Airspace is then assigned to each flight cluster. An improved version of \textit{SectorFlow} lets the algorithm keep flow intersections away from sector boundaries and use a gradient search to refine the boundaries after the initial partition\cite{brinton2009airspace}. Meanwhile, \textit{CellGeoSect\cite{sabhnani2010flow}} is a cell clustering method that represents the airspace as a group of hexagonal cells. Then clusters these cells to maximize flow connectivity and balance flight count between each cluster. It then modifies the resulting design by removing and redefining the boundary between each pair of airports to avoid the geometric constraints of major flows. It is important to note that air traffic controllers need the training to work on a specific airspace set\cite{zelinski2011comparing}. It is not desirable operationally since ATC personnel cannot familiarize themselves with the newly assigned airports\cite{Sergeevagenetic}. 

First, dynamically creating entirely new controlled airspace. Controlled airspace sector boundaries are newly created every time without relying on pre-existing structures. This approach is also referred to as dynamic sectorization. It is important to note that air traffic controllers need the training to work on a specific airspace set\cite{zelinski2011comparing}. So that previous work mentioned\cite{brinton2008airspace, sabhnani2010flow} following this path is not desirable operationally since ATC personnel cannot familiarize themselves with the newly assigned airports\cite{Sergeevagenetic}. 

The second approach is to stick to existing building blocks, e.g., airspace modules, that can be dynamically combined to form a controlled airspace sector\cite{Sergeevagenetic}. Currently, this approach is more desirable operationally. A controlled airspace sector is operated by a small team of controllers and comprised of one or more airspace modules. In both the U.S. and European airspace, controlled airspace can be combined with others or split into smaller controlled regions to balance the workload equally across available ATC resources and airports\cite{zelinski2011comparing}. Several examples incorporating various constraints are as follows:

% While limited research satisfies the discussed factors, some research frameworks offer promising approaches for realizing and optimizing DAC. These frameworks utilize fixed airspace modules that provide a stable foundation for air traffic controllers to work with, minimizing safety risks associated with frequent airspace configuration changes. 

% tree search
In \cite{gianazza2010}, a promising framework is proposed to realize DAC. First, features for evaluating air traffic controllers' workload are extracted from the flight radar track, and sector operation history. These feature vectors are fed into a neural network to provide a workload indication for the ATC in terms of high, normal, or low. Their algorithm generates different new configurations by splitting airports into several smaller airspace modules when the workload is high or merging with other airports when the workload is low. Next, the tree search methods explore all possible partitions while restricting them to be operationally valid. This ensures the algorithm builds an optimal airspace partition where the workload is balanced across the airports and uses the restrictions to lower the reconfiguration cost. 

% genetic algorithm
Sergeeva et al.\cite{Sergeevagenetic} proposed to model dynamic airspace configuration as a graph partitioning problem that can be optimized with a genetic algorithm. They define two different types of airspace modules. Those airspace modules that “are permanently busy areas with a high traffic load” are designated to be “Sector Building Blocks” (SBBs). Less busy and more generic airspace modules are called “Sharable Airspace Modules” (SAMs). A controlled airspace sector should consist of at least one SBB and multiple SAMs. Instead of incorporating re-configuration into the cost function, this approach ensures the stability of the configuration by making the busiest airspace modules (SBBs) a fixed central component of each controlled airspace sector. Only the generic SAMs change from one configuration to the next. The approach works best when the airspace is divided into small SAMs and SBBs. 

Spectral clustering is applied in \cite{li2009spectral} and \cite{Sergeevagenetic} to balance air traffic controllers' workload. These works first transform the airspace configuration problem into a graph partitioning problem and address it with spectral clustering. The main idea of spectral clustering is first to perform eigendecomposition on the adjacency matrix to extract essential components, then use the k-means clustering algorithm to divide the airspace. The graph uses vertices and links to model airports, waypoints, and air routes, then project real flight trajectories onto the graph as edge weight. The key idea of spectral clustering is to reduce the graph dimension to help the clustering algorithm focus on the most critical feature. However, eigendecomposition is a linear operation that may not have the flexibility to capture the main component of complex air traffic patterns \cite{gondara2016medical}. 

In general, the previous works contains the following drawbacks: (a) Omiting the consideration that only close airport should collaborate. (b) Involving complicated hyperparameters that are not self-adjusted. (c) The time complexity is high. This is the primary motivation for our research.  
% The previous works have the following drawbacks:
% \begin{itemize}
%     \item The clustering algorithms require the configuration of parameters that do not correlate to practical metrics.
%     \item The overhead for dimension reduction is large.
% \end{itemize}

%Compared with eigendecomposition, autoencoder is a nonlinear model used a lot to compress and denoising images\cite{vincent2008extracting} and has shown a strong ability to extract the most effective features, which can be used to generate new images\cite{mescheder2017adversarial}

\section{Methodology}
\label{sectMM}
\subsection{Problem Definition \& Datasets}
We try to find an optimal plan to allocate non-busy airports’ air traffic control resources to assist busy airports with possibly more delays. Our goal is to balance the workload of different airports during emergency evacuations or other busy scenarios. We use a metric called the Regional Unbalanced Level (RUL) to quantify the workload of handling regular or delayed flights. We first calculated the average number of non-delayed and delayed flights handled by airports within each cluster, where m is the total number of flights in the cluster \(i\):
\begin{align}
    F_{i} = \frac{\sum_{k=1}^{m}f_{k}}{m}\\
    D_{i} = \frac{\sum_{k=1}^{m}d_{k}}{m}
\end{align}
    
We then calculated the variance of the \(F_{i}\) and \(D_{i}\) over all clusters, noted as \(S\). This value helps to quantify the workload unbalance within the whole airspace in scope. \(n\) is the number of clusters in the configuration, D and F are the mean workload of delayed and on-time flights: 
\begin{align}
    S_{D} = \frac{\sum_{i=1}^{n}(D_{i}-D)^2}{n-1}\\
    S_{F} = \frac{\sum_{i=1}^{n}(F_{i}-F)^2}{n-1}
\end{align}
Where high \(S_{D}\) or \(S_{F}\) indicate the workload of handling delayed or regular flights are highly different across different clusters. We assume that airports within the same cluster are collaborating with each other to handle emerging workloads. Therefore, smaller \(S_{D}\) or \(S_{F}\) are preferred.

%We believe this method can be generalized not only to emergency situations but can also be used in daily routines when passenger flows are distributed in a temporally and spatially unbalanced manner.

\subsection{Dataset}
Flight Delays and Cancellations were published by The U.S. Bureau of Transportation in 2015\cite{Transportation_2017}\textcolor{blue}. This dataset includes statistics tracking the on-time performance of domestic flights operated by large air carriers. The original data has a total of 30 attributes, however, not every attribute was recorded for each flight, thus columns with more than 25\% missing values are removed. For the rest of the data, only the related attributes are kept, including the scheduled date, airline, origin and destination airport, departure time, and delay time. Canceled or diverted flights are removed. To reduce computation load, departure delay is transformed into the binary label, 0 is a non-delay flight and 1 is a delayed flight. We focused on the 21 airports in Florida. 

\subsection{Hybrid Graph Modeling for Airspace}
\begin{figure}[h!]
    \centering
    \includegraphics[width = 0.7\linewidth]{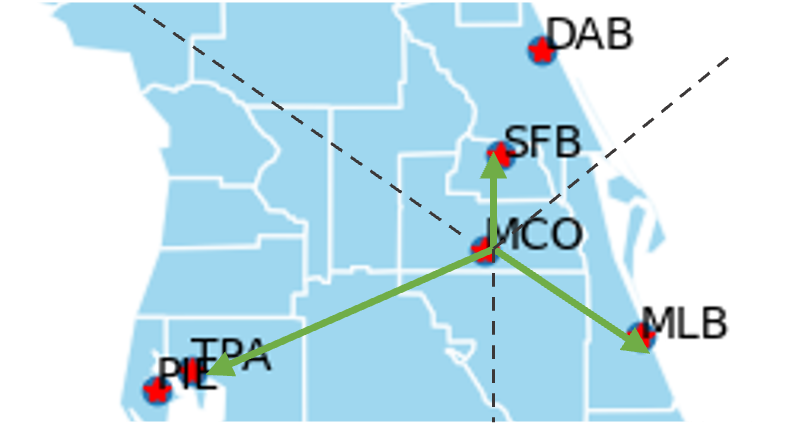}
    \caption{Illustration of Initial Airport Adjacency Graph (IAG). The example airports are from central Florida.}
    \label{figIAG}
\end{figure}

Modeling air traffic system in the form of a graph can effectively preserve the spatial and temporal information in the system\cite{jiang2021spatial}. We used a novel data structure Hybrid Airport Adjacency Graph (HAG) to model the airspace of incorporating geographical adjacency and workload-based mergeability. The procedures are as follows:

\textbf{Step 1:} \textit{Generating Initial Airport Adjacency Graph (IAG)}:
We identify if two nodes are connectible based on the geographic location. Mathematically, we define that if two airports, $V1$ and $V2$ are connected if $V2$’s is the closest neighbor geographical neighbor of $V1$ at the same azimuth. Figure \ref{figIAG} is an example taken from central Florida. When setting the azimuth to 120 degrees, node MCO has three connected nodes which are SFB, TPA, and MLB. Whereas, PIE, DAB are not connected to MCO because they are not the closest node inside the 120 azimuths. The size of the azimuth is a tuning parameter, where a smaller azimuth resulting a more dynamic graph, but may result in connecting to a non-realistic node that is too far away. 

\textbf{Step 2}: \textit{Creating Hybrid Airport Adjacency Graph (HAG):} 
We assigned edge weights for each connection of airports. In HAG, airports that are less busy should have stronger connections to each other, and the busy airports should be far away from each other to avoid being clustered together. In other words, the edge weight is negatively related to the total workload quantified by the estimated delay flights and delay ratio in the future two hours between the nodes. A modified radial-based kernel is chosen to encode the workload into edge weight, defined as below:
\begin{align}
    load_{ij} &= \frac{d_{i}+d_{j}}{f_{i}+f_{j}} \\
    w_{i,j}  &= B^{((1-\lambda)(100load_{ij})+\lambda d_{ij}-shift)}
\end{align}
where $i$ and $j$ are two connected nodes, $d$ is the number of delayed flights in this time window, $f$ is the total number of flights in this time window. Thus, $load_{ij}$ is the percentage of delayed flights of node $i$ and $j$ and normalized between 0 and 1. We only use the percentage of delayed flights as an indicator of the ATC workload, a more comprehensive workload can be explored in the future and substituted here. $d_{ij}$ is the geographical distance between the two airports. $\lambda$ is the geographical weight factor in balancing between considering geographical distance and gross delay ratio, this factor is automatically adjusted in our program. $B\in (0,1)$ is the base, such a base is to satisfy the negative relation between workload and the edge weight. The sift is set as 300 as it is numerically more stable. Mathematically, the larger the edge weight between two nodes (airports), the more likely they should be connected to collaborate and form a cluster.

\subsection{Adaptive Clustering for Pre-Allocation}
We developed an adaptive clustering algorithm for partitioning the fine-tuning of the collaboration airports, with Spectral Clustering\cite{ng2001spectral} as a key component. Mathematically, spectral clustering first performs Eigen Decomposition of the adjacency matrix of the Hybrid Airport Adjacent Graph to project data from a higher dimension to a lower dimension to remove redundancy and noise, then, clustering is done on the low-dimension representation of data \cite{ng2001spectral}. Our procedures are as follows:

\textbf{Step 1}: Construct Hybrid Airport Adjacency Graph with \(\lambda = 0\), in this way, the initial clustering will not consider the geographical location of airports.

\textbf{Step 2}: Compress the HAG's adjacency matrix, we explore the following methods:
\begin{itemize}
    \item SVD: we calculate the degree matrix of the graph; The degree matrix is a diagonal matrix where the value at entry \((i,i)\) is the degree of the node \(i\). Then calculate the eigenvalues and eigenvectors of the degree matrix; then we sort them based on the eigenvalues. 
    \item Autoencoder: we use the encoder to compress each airport in the graph (each row) into lower dimensions. This autoencoder contains two Dense layers for its encoder and decoder respectively.
\end{itemize}

\textbf{Step 3}: We perform the k-means clustering algorithm with an initial $k$ value equal to half of the airports on the low-dimension graph to get the initial clustering result.

\textbf{Step 4}: Scan each cluster in the initial clustering result, if any cluster contains more than three airports or with a diameter greater than 100 nautical miles (the typical transmission range that aircraft can communicate directly). We increase the number of clusters by 1 and simultaneously, increase the geographical weight  $\lambda$ by 0.1, but $\lambda$ can not exceed 0.5.

\textbf{Step 5}: We repeat the clustering procedure as described in Step (1) until all clusters satisfy the criteria defined in Step 4.

Steps 4 and 5 make the clustering process adaptive and this algorithm can gradually evolve to use geographical constraints to create clusters with reasonable spatial size. Therefore, we do not require select dedicate  \(k\) and  \(\lambda\) values for different scenarios.

After spectral clustering on HAG, the airports that are geographically close and with relatively low workloads are combined as a new cluster. Simultaneously, the busy airports with more delayed flights will be picked up and isolated. 

\subsection{Fine-Tuning for Dynamic Workload Balancing}
In this stage, we aim to merge different airports’ governing regions to rebalance the workload of the area when there is an abrupt increase in travel demand or flight delays. For this purpose, we have the following assumptions:

\textit{\textbf{Assumption I: }{All flight plans are known at least two hours in advance from flight plans or predictions.}}

\textit{\textbf{Assumption II: } The abrupt increase in travel demand under emergency situations could cause significant delays in flights. In our experiment, if an airport's delayed flights within the predicted time window are greater than 2 delayed flights per hour (with a regional airport) or with a delay percentage within this time window being 30\% (medium or large airports), we then mark this airport as a busy airport and needs external assistance.}

\textit{\textbf{Assumption III: } The nearby airport that used to assist a busy airport should: (a) have fewer delayed flights if the category of the airports is identical, or a lower percentage of delay if the category of the airport is different.}

Based on assumptions II and III, we develop the fine-tuning algorithm for each busy airport as follows:

\begin{itemize}
    \item \textbf{Step 1}: we created a ranked list of busy airports based on (a) a user-defined priority level with a default value of zero, (b) the number of delayed flights, (c) delay ratio, and (d) number of scheduled flights within the predicted time window. Here, the user-defined priority level can be filled when there’s an emergent situation.
    \item \textbf{Step 2}: we scanned all airports within 100 nautical miles of the busy airports in step one and determine if a specific airport can be merged to assist an adjacent busy airport based on these criteria: (a) distance, (b) less number of predicted delayed flights, (c) lower delay ratio. Specifically, we created a ranked list of non-busy candidates and picked the closest one.
    \item \textbf{Step 3}: if any two airports are selected as a collaborative pair, we created a new cluster containing only the two airports, to prevent airspace conflict, we also ensure that there’s no other busy airport within the combined airspace before establishing the collaboration relationship.
\end{itemize}
In general, the algorithm allows reallocating more resources from less busy regions in the airspace.

% \begin{figure}[b]
%     \centering
%     \includegraphics[width = 0.9\linewidth]{Figures/1224_diff.png}
%     \caption{New configuration generated at different times on 12/24 (a) 7:00-9:00, (b) 19:00-21:00 }
%     \label{fig1224_diff}
% \end{figure}

% \begin{figure}[hbt!]
%     \centering
%     \includegraphics[width = 0.5\textwidth]{Figures/1224_diff.png}
%     \caption{New configuration generated at different times on 12/24 (a) 7:00-9:00, (b) 12:00-14:00, (c) 19:00-21:00 }
%     \label{fig1224_diff}
% \end{figure}

\begin{figure}[hbt!]
\centering  
\subfloat[]
{%
    \includegraphics[width=0.7\linewidth]{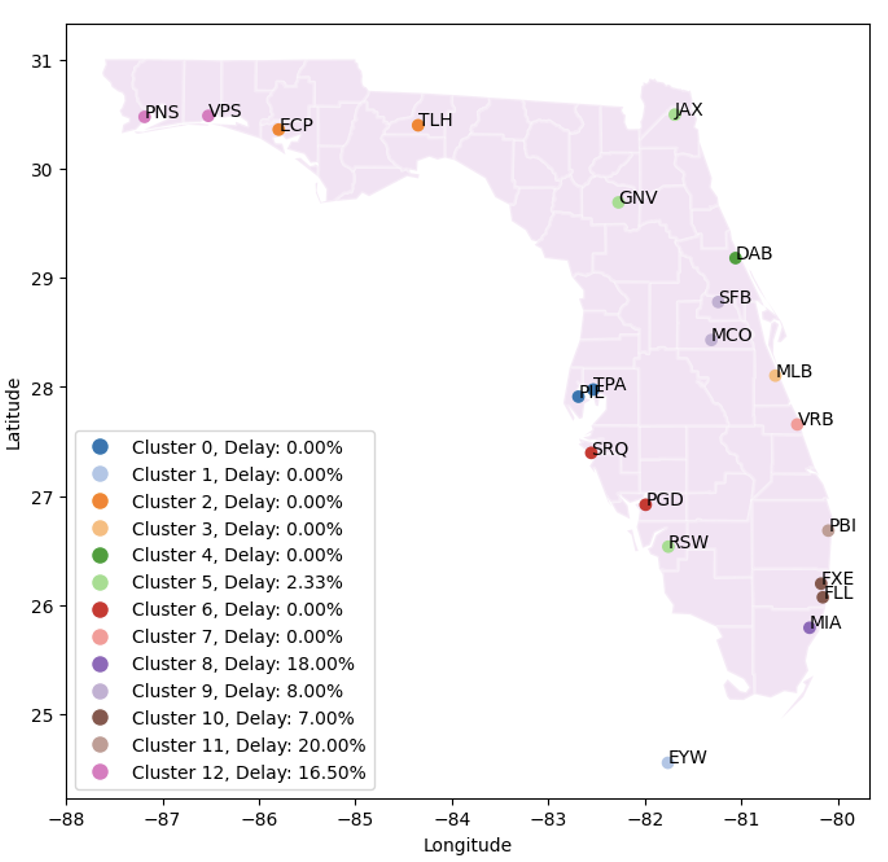}

}\\
\vspace{-10px}
\subfloat[]
{
    \includegraphics[width=0.7\linewidth]{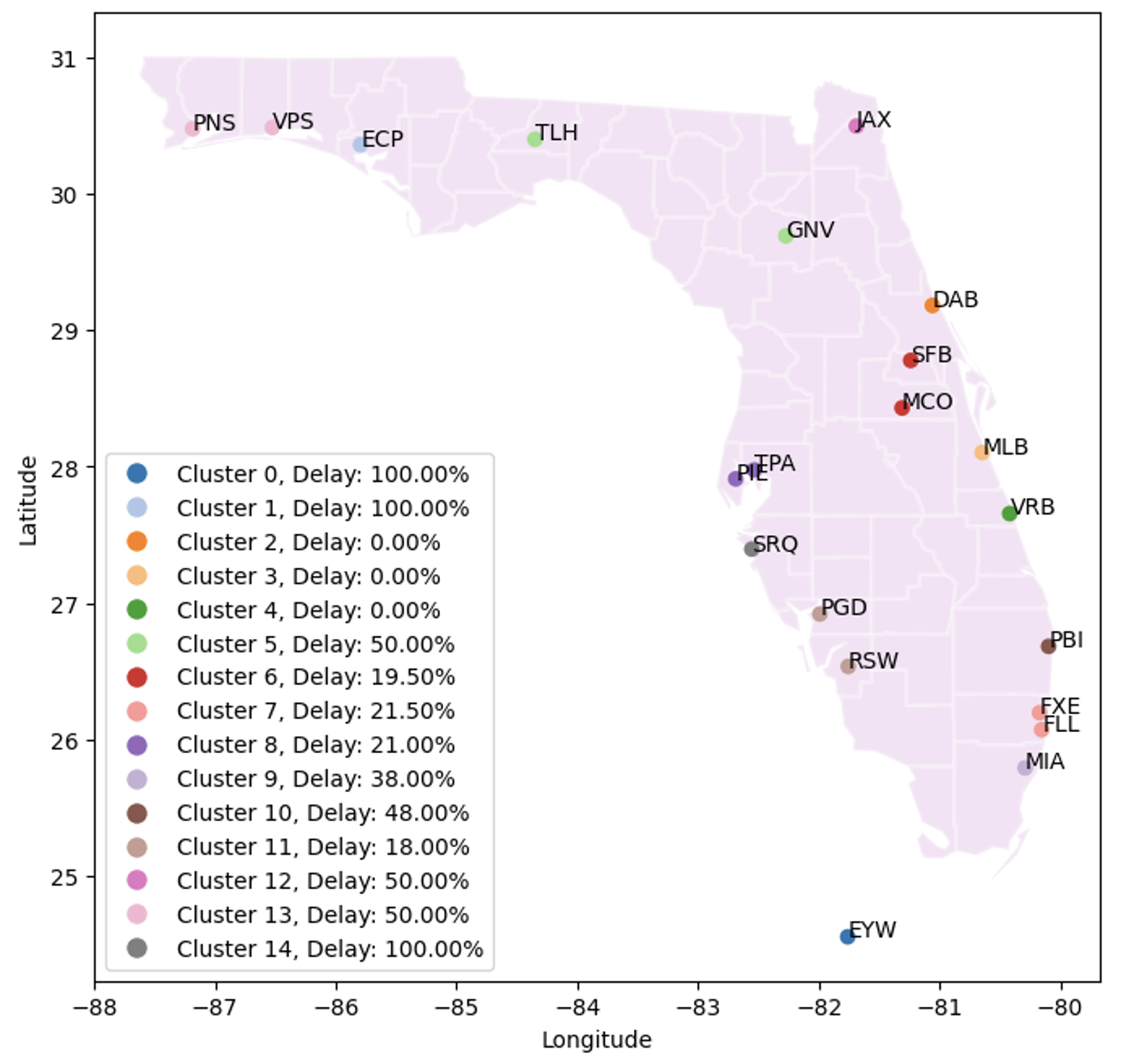}
}%
\caption{Airspace configuration generated at different times on 12/24 (a) 7:00-9:00, (b) 12:00-14:00.}
\label{fig1224_diff}
\end{figure}
%%%%%%%%%%%%%%%%%%%%%%%%%%%%%%%%%%%%%%%%%%%%%%%%%%%%%%%%%%%%%%%%%%%%%%%%%%%%%%%%
\section{Evaluation \& Discussion}
\label{sectEED}
This section examines the efficiency of the graph clustering-based dynamic airspace configuration method under different scenarios. Additionally, we compare the computational efficiency of Autoencoder and SVD for graph compression.

\subsection{DAC with different scenarios}

\subsubsection{Different times on the same day}

We tested the DAC algorithm on December 24th, 2015, a day with heavy flight traffic due to the holiday season. Three distinct 2-hour time windows are selected: 7:00-9:00 for low traffic, 12:00-14:00 for high traffic, and 19:00-21:00 for medium traffic. The experiment results show that the algorithm changes the configuration based on the different traffic conditions successfully, and the workload unbalanced level of the new configuration in terms of handling regular and delayed flights has been significantly reduced at all different traffic conditions, as seen in Table \ref{table1224_diff}.

\begin{table}[h]
\caption{Reduction on unbalance level after reconfiguration at different times on the same day}
\label{table1224_diff}
\begin{adjustbox}{width=\columnwidth,center}
\centering
\begin{threeparttable}

\begin{tabular}{@{} l *{2}{C{2.5cm}} @{}}
\toprule
&  Handling regular flights
 &  Handling delayed flights\\
\midrule
7:00-9:00 (Low traffic) & 42.85\% & 10.8\% \\ 
12:00-14:00 (High traffic) & 56.9\% & 61.04\% \\ 
19:00-21:00 (Medium traffic) & 42.86\% & 60.1\% \\
\bottomrule
\end{tabular}
\end{threeparttable}
\end{adjustbox}
\end{table}

As in Figure~\ref{fig1224_diff}, for low traffic conditions between 7:00-9:00, the airspace configuration algorithm combines adjacent sectors to balance the workload, such as MCO with SFB, FXE with FLL, and PNS with VPS; In the meantime, several airports are isolated independently because there is no non-busy airport within a reasonable range and without airspace overlap with busy airports, such as MIA and EYW. When the traffic load is high during 12:00-14:00, the algorithm changes the configuration to balance the workload among sectors to adapt to the increasing traffic load. 

% \begin{figure*}[h]
%     \centering
%     \includegraphics[width = 0.85\textwidth]{Figures/diffday_1214.png}
%     \caption{New configuration generated on different days for (a)high and (b)low traffic load (a1)7/3, (a2)11/25(one day before Thanksgiving), (a3)12/31; (b1)2/17, (b2)6/9, (b3)9/8}
%     \label{figdiffday}
% \end{figure*}

\begin{figure}[h]
\centering
\hspace{-20px}\subfloat[]
{%
    \includegraphics[width=0.8\linewidth]{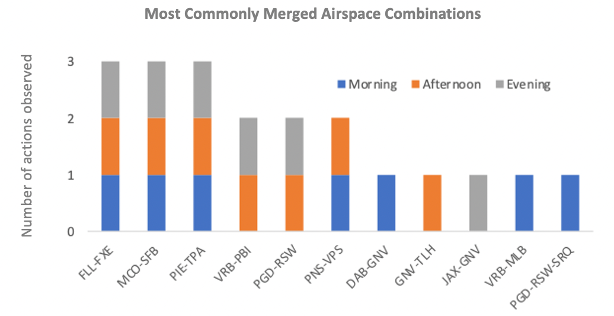}
}\\
\vspace{-10px}
\hspace{-20px}\subfloat[]
{
    \includegraphics[width=0.8\linewidth]{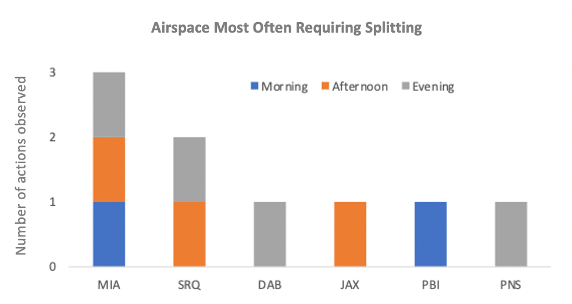}
}%
\caption{Summary of airspace reconfiguration actions within Dec 24th 2015  (a) Merged airspace. (b) Airspaces that need further separation.}
\label{fig1224_analysis}
\end{figure}

We also noticed some common patterns. First, the airspace merging actions in the three time periods are summarized in Figure~\ref{fig1224_analysis}(a). The three most frequent merges are FLL-FXE, MCO-SFB, and TPA-PIE. These airports are selected to merge into collaborative pairs because a) executive airports are not usually as busy as large international airports even if they are closely located, so they can always assist busy airports. Secondly, Figure \ref{fig1224_analysis}(b) shows that MIA and its nearby airports are extremely busy all day round, which makes it impossible to find a collaborative airport that is less busy. 

\begin{figure}[hbt!]
	\centering  
	\subfloat[]
	{%
		\includegraphics[width=0.8\linewidth]{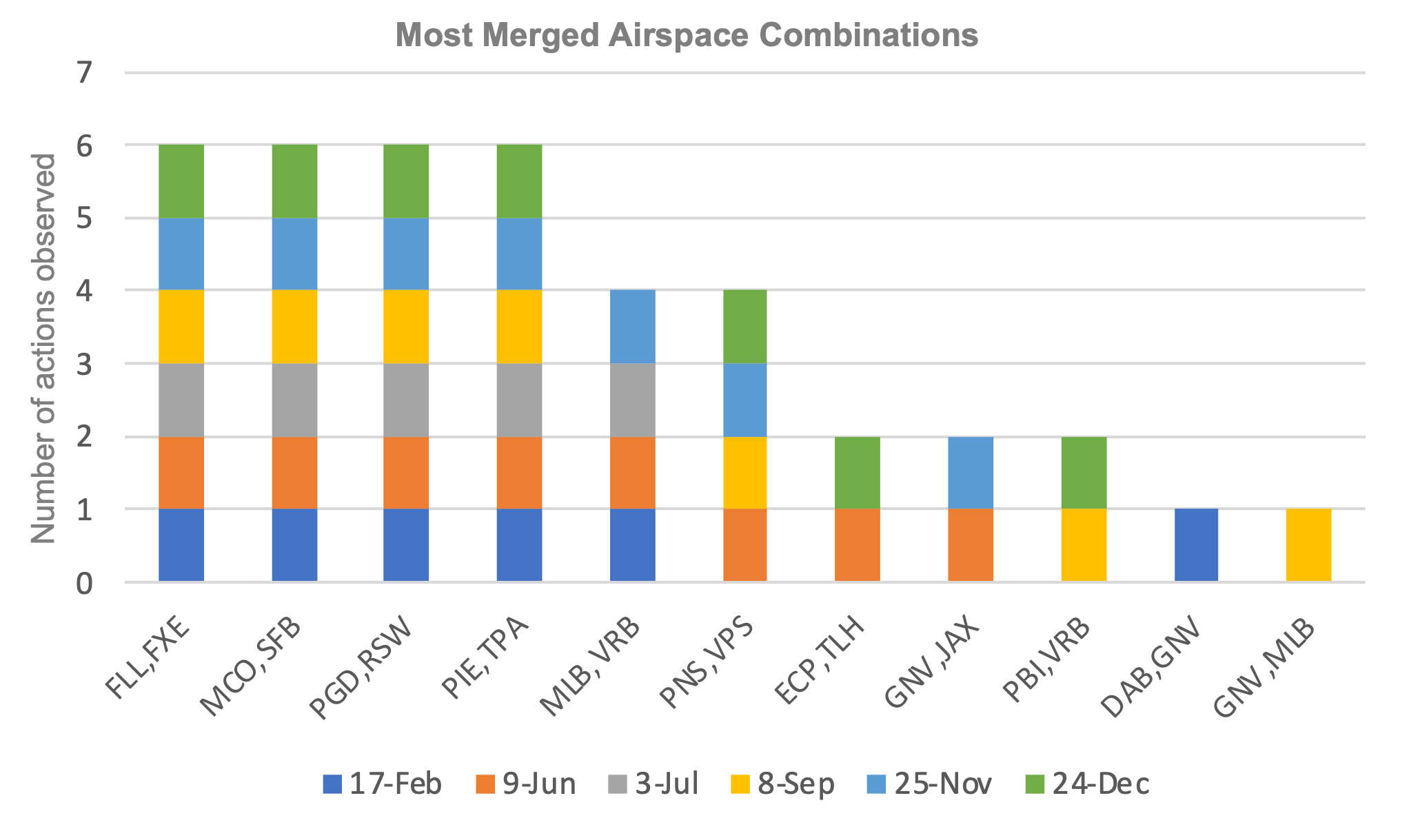}
	}\\
	\vspace{-10px}
	\subfloat[]
	{
		\includegraphics[width=0.8\linewidth]{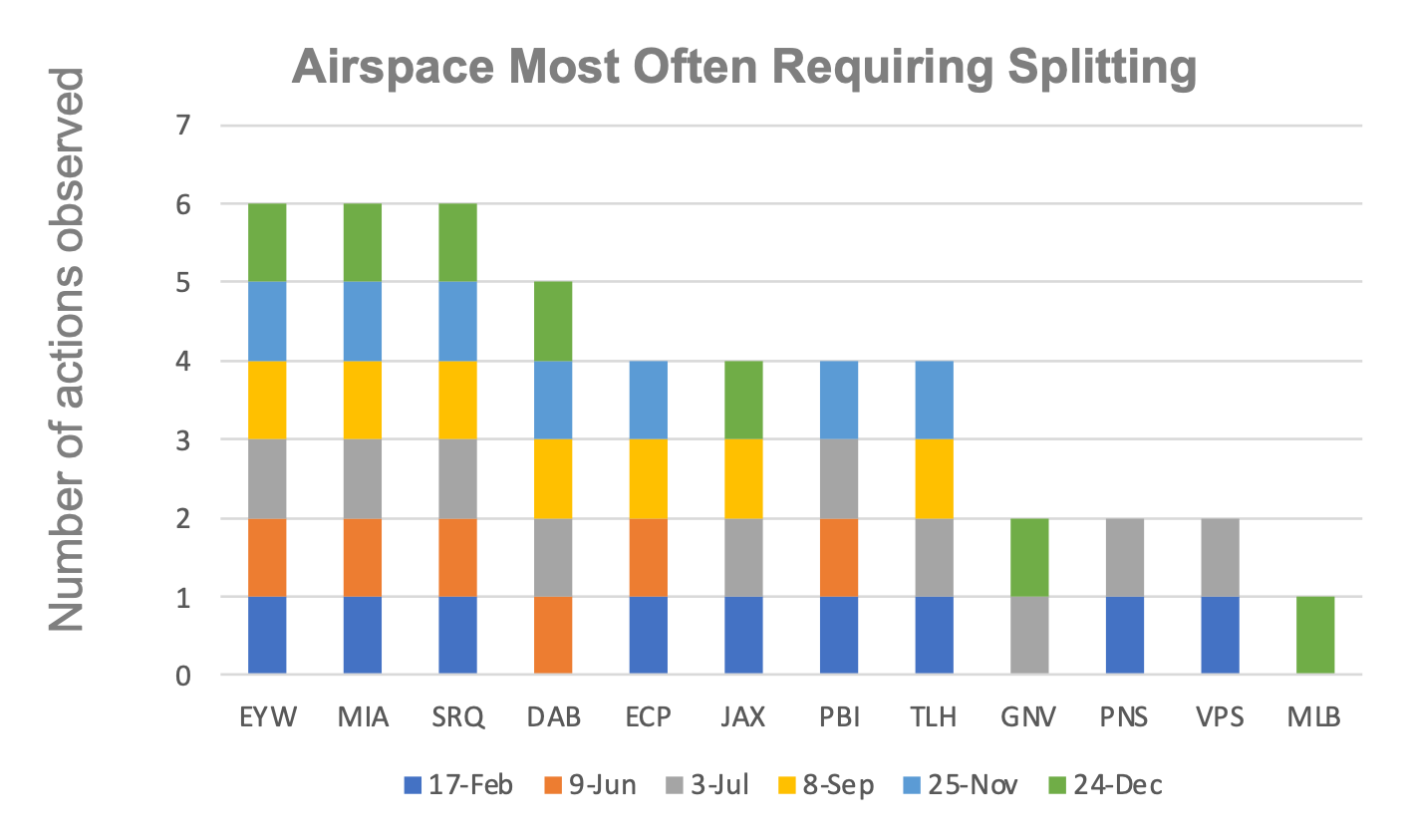}
	}%
	\caption{Summary of airspace reconfiguration actions observed on 12:00 PM to 14:00 PM on different days of 2015: (a) airspaces that are merged. (b) airspaces that need further separation.}
	\label{figdiffday_analysis}
\end{figure}

\begin{figure*}[b!]
	\centering
	\includegraphics[width = 0.95\textwidth]{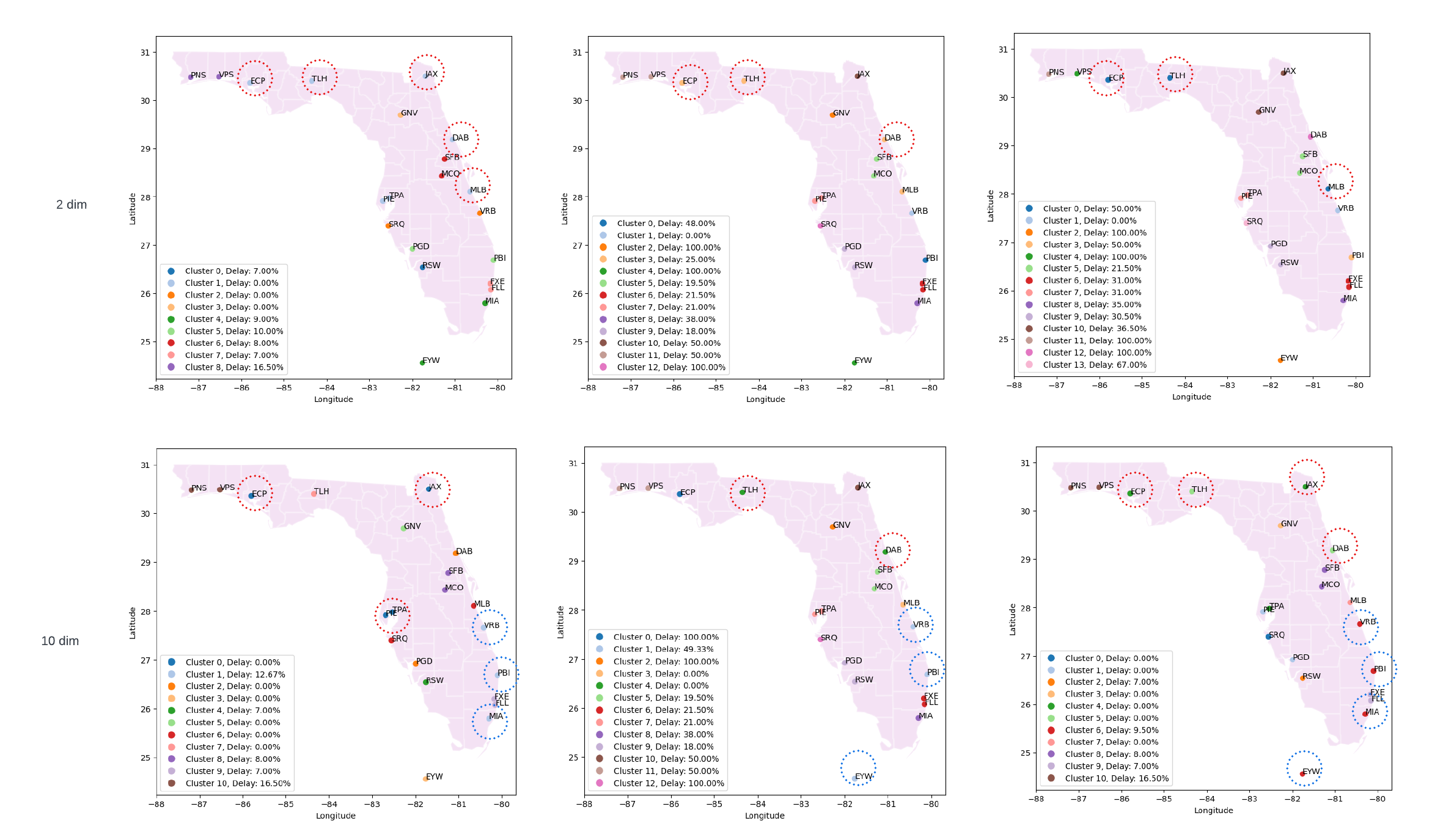}
	\caption{Error patterns of autoencoder in pre-clustering graph data compression. The result of the first row and second row are generated by AE when the latent space is 2 dimensions and 10 dimensions respectively.}
	\label{figAEFlawpattern}
\end{figure*}

\subsubsection{Same time for different days}
We also evaluate the algorithm on different dates, three high-traffic volume dates are selected, which are 7/3, 11/25 (one day before Thanksgiving), and 12/31 (one day before the new year) in 2015. we also compare the airspace configuration results on the busy hours (12:00 PM to 14:00 PM) on these low-traffic days: 2/17, 6/9, and 9/8, 2015. Our algorithm significantly decreases unbalance level of the ATC’s workload in both busy and non-busy scenarios as shown in table \ref{tablediffday}.

\begin{table}[hbt!]
\begin{adjustbox}{width=\columnwidth,center}
\centering
\begin{threeparttable}
\caption{Reduction on unbalance level after reconfiguration at the same times on different days}
\label{tablediffday}

\medskip

\begin{tabular}{@{} l *{2}{C{2.5cm}} @{}}
\toprule
High traffic dates & Handling regular flights
 &  Handling delayed flights \\
\midrule
7/3 & 56.77\% & 56.37\% \\ 
11/25 & 59.1\% & 62.5\% \\ 
12/31 & 56.9\% & 59.9\% \\
\bottomrule
\end{tabular}

\bigskip

\begin{tabular}{@{} l *{2}{C{2.5cm}} @{}}
\toprule
Low traffic dates & Handling regular flights
 &  Handling delayed flights \\
\midrule
2/17 & 64.1\% & 66.7\% \\ 
11/25 & 60\% & 67.1\% \\ 
12/31 & 53\% & 44.5\% \\
\bottomrule
\end{tabular}
\end{threeparttable}
\end{adjustbox}
\end{table}

There are also some common merging strategies employed for both high and low traffic conditions. As shown in figure \ref{figdiffday_analysis}, the actions taken most frequently are the merging between FLL-FXE, MCO-SFB, PGD-RSW and PIE-TPA, showing that large international airports are always busy all year round and need assistance from other regional or executive airports. GNV is surrounded by medium and busy airports and is dynamically assigned as a flexible collaborator. Our algorithm also indicates that some busy airports in southern Florida are surrounded by airports that are busy simultaneously. Consequently, they can not find collaborative peers and further divisions are needed in their internal airspace.

\subsection{Comparison of pre-clustering graph data compression methods}
\label{sectEvaBehvDev}

\begin{figure}[hbt!]
\centering  
\subfloat[]
{%
    \includegraphics[width=0.7\linewidth]{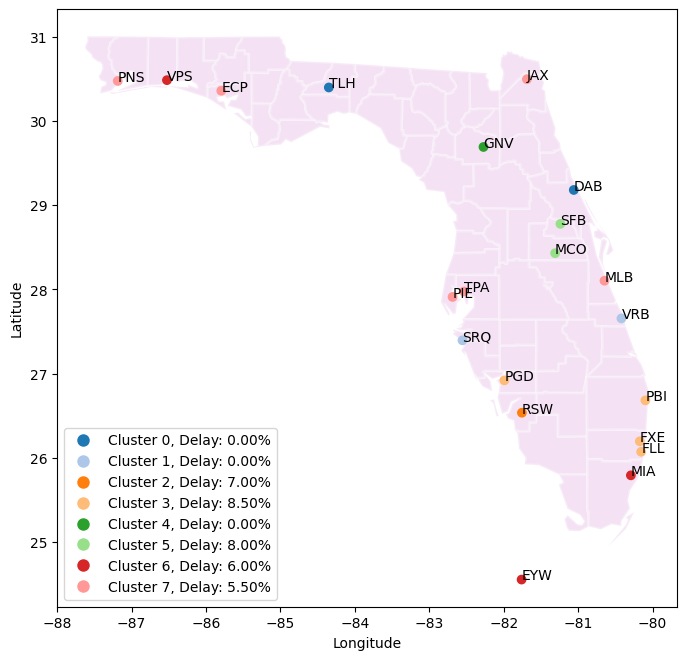}
}\\
\vspace{-10px}
\subfloat[]
{
    \includegraphics[width=0.7\linewidth]{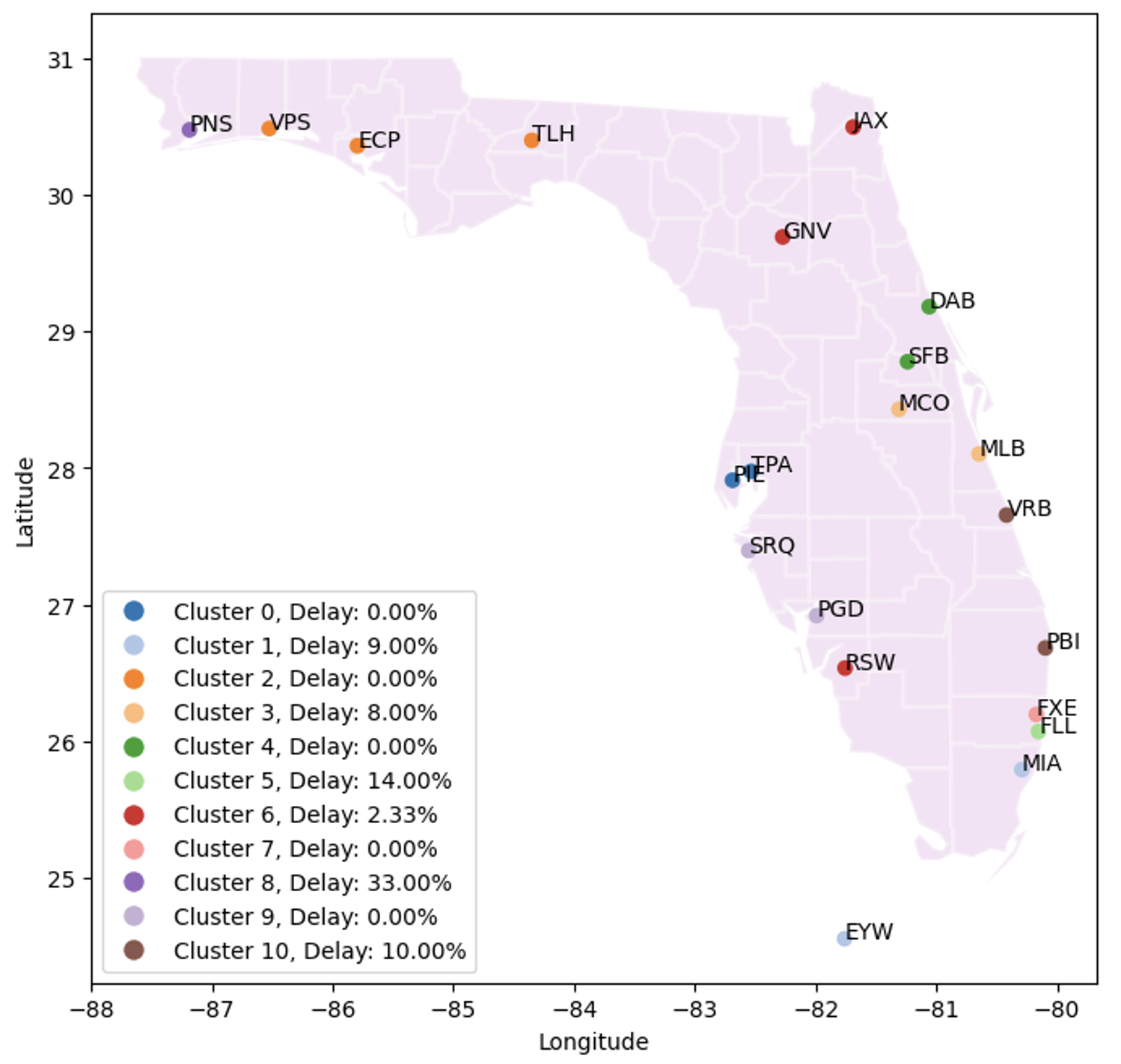}
}%
\caption{Comparison of pre-clustering compression method on the final result of DAC (a) Nonlinear: Autoencoder and (b) Linear: SVD.}
\label{figAEvsSVD}
\end{figure}

SVD and Autoencoder (AE) are employed to reduce the dimension of the HAG's adjacency matrix. In general, both techniques yield low-dimension graphs that preserve a significant amount of essential information, but each has its own advantages. Specifically, SVD excels in preserving the embedded constraints within the graph, whereas the AE is particularly adept at achieving balanced clustering. In Figure \ref{figAEvsSVD}, the clustering results of the SVD-compressed graph adhere to the adjacency constraint that we only want to merge geographically close airports. Comparably, an AE allows the algorithm to deviate slightly. For instance, when using the AE depicted in (a) for compression, certain clusters, such as clusters 0 and 7 contain airports that are not directly connected. 

To compare the impact of the dimension of the latent space to the strictness of AE, multiple AEs were trained with varying latent dimensions (2, 5, 10 and 15). We found that there is a certain pattern when AEs do not follow the restrictions. As in figure \ref{figAEFlawpattern}, AEs always mistakenly combine ECP, TLH and DAB together or VRB, PBI, MIA, and EYW together, but this kind of error can easily be fixed by a post-processing algorithm. %Furthermore, when assessing the level of workload imbalance, the autoencoder significantly outperforms SVD in achieving balanced workload distribution among clusters. In fact, the clustering result using the autoencoder is approximately seven times less unbalanced than the result obtained through SVD.

% \begin{figure}[hbt!]
%     \centering
%     \includegraphics[width = 0.7\linewidth]{Figures/AE_flawpattern.pdf}
%     \caption{Error patterns of autoencoder in pre-clustering graph data compression. The result of the first row and second row are generated by AE when the latent space is 2 dimensions and 10 dimensions respectively.}
%     \label{figAEFlawpattern}
% \end{figure}

Additionally, AEs are much faster than SVD after training, making them more suitable for use in real time. As in table \ref{tableAESVD_time}, as the dimension of the graph continues to increase, SVD is much slower than AEs. This is because the computation of trained AEs can easily get accelerated by hardware. Also, AEs can leverage past training experiences, unlike SVD which starts computation from scratch each time.

This is because the computation of trained AEs can easily get accelerated by hardware. Also, SVD has to start computation from scratch every time while AEs can leverage past experiences from training. 

\begin{table}[hbt!]
\caption{Comparison of time required in real-time application between SVD and Autoencoder }
\label{tableAESVD_time}
\begin{adjustbox}{width=\columnwidth,center}
\centering
\begin{threeparttable}

\begin{tabular}{@{} l *{4}{C{2cm}} @{}}
\toprule
 HAG's dimension & SVD prediction ($10^{-5}$ second) &  Autoencoder prediction ($10^{-5}$ second) & Autoencoder Training time (second)\\
\midrule
7*7 & 43.9 & 1.16 & 27.78\\ 
10*10 & 44.8 & 1.70 & 27.04\\ 
15*15 & 44.5 & 4.06 & 69.47\\
18*18 & 48.1 & 3.30  & 64.67\\
21*21 & 373.8 & 3.54 & 65.98\\
\bottomrule
\end{tabular}
\end{threeparttable}
\end{adjustbox}
\end{table}

% \begin{tabular}{@{} l *{2}{C{2.5cm}} @{}}
% \toprule
% Low traffic dates & Handling regular flights
%  &  Handling delayed flights \\
% \midrule
% 2/17 & 64.1\% & 66.7\% \\ 
% 11/25 & 60\% & 67.1\% \\ 
% 12/31 & 53\% & 44.5\% \\
% \bottomrule
% \end{tabular}
% \end{threeparttable}
% \end{adjustbox}
% \end{table}

% \subsection{Limitations of this research}

% There are several limitations to our study at the current stage: (a) Since we were not able to fetch detailed information on ATC controllers in many airports, it makes it impossible to create a precise model to quantify the capability of each tower and balance their workloads. (b) We leveraged the delay ratio as a metric to roughly estimate the busy level of airports, better metrics can be used to make the model more practical. (c) The validity of data from public sources is quite limited, as airports’ operational data are seldom archived and shared properly, which, in turn, makes the simulation less realistic and applicable.

%%%%%%%%%%%%%%%%%%%%%%%%%%%%%%%%%%%%%%%%%%%%%%%%%%%%%%%%%%%%%%%%%%%%%%%%%%%%%%%%%%%%%%
\section{Conclusion}
We've introduced an innovative method for dynamic airspace reconfiguration, using a graph model to balance air traffic controllers' workload. This model incorporates geographical adjacency and ATC workloads, utilizing a spectral clustering-enabled adaptive algorithm to generate new configurations based on predicted delay and flight plans. The algorithm groups high-workload airports as centers surrounded by lesser-engaged airports, redistributing ATC resources for workload equilibrium. The model outperforms static airspace configurations, reducing workload imbalance by over 55\% under high traffic volume, as confirmed through evaluations during various time windows and traffic conditions. Our real-data simulations indicated that Miami and Sarasota airspace requires further partitioning to improve performance, as our DAC algorithm struggled to locate less busy airports within 60 miles for collaboration. Furthermore, Key West airport's remote location hindered collaboration with busy airports to share emergent workloads.

% Our simulations on real data also revealed that the air space in Miami and Sarasota needs to be further sliced to increase its capability since our DAC algorithm has a low success rate in finding collaboratively airports that are less busy within 60 miles.  Other than that, the airspace in Key West airport locates too far away from other airports making it almost impossible to collaborate with other busy airports to share emergent workloads. 

% There are several limitations to our study at the current stage: (a) Since we were not able to fetch detailed information on ATC controllers in many airports, it makes it impossible to create a precise model to quantify the capability of each tower and balance their workloads. (b) We leveraged the delay ratio as a metric to roughly estimate the busy level of airports, better metrics can be used to make the model more practical. (c) The validity of data from public sources is quite limited, as airports’ operational data are seldom archived and shared properly, which, in turn, makes the simulation less realistic and applicable.

Our future direction includes improving airspace ATC workload assessment by considering metrics beyond delayed and total flights. Additionally, we plan to explore neural networks' potential in generating comprehensive airspace configuration plans.

\label{sectCC}

\section*{Acknowledgment}

This research was supported by the Center for Advanced Transportation Mobility (CATM), USDOT Grant No. 69A3551747125, 270128BB(AWD00237), and by the U.S. National Science Foundation under Grant No.2231629, Grant No.2142154, No.2142514 and Grant No.2309760.

\bibliographystyle{IEEEtran}
% Tell the complier which style you want.
\bibliography{IPCCC2022.bib}

% Generated by IEEEtran.bst, version: 1.14 (2015/08/26)
\begin{thebibliography}{10}
\providecommand{\url}[1]{#1}
\csname url@samestyle\endcsname
\providecommand{\newblock}{\relax}
\providecommand{\bibinfo}[2]{#2}
\providecommand{\BIBentrySTDinterwordspacing}{\spaceskip=0pt\relax}
\providecommand{\BIBentryALTinterwordstretchfactor}{4}
\providecommand{\BIBentryALTinterwordspacing}{\spaceskip=\fontdimen2\font plus
\BIBentryALTinterwordstretchfactor\fontdimen3\font minus
  \fontdimen4\font\relax}
\providecommand{\BIBforeignlanguage}[2]{{%
\expandafter\ifx\csname l@#1\endcsname\relax
\typeout{** WARNING: IEEEtran.bst: No hyphenation pattern has been}%
\typeout{** loaded for the language `#1'. Using the pattern for}%
\typeout{** the default language instead.}%
\else
\language=\csname l@#1\endcsname
\fi
#2}}
\providecommand{\BIBdecl}{\relax}
\BIBdecl

\bibitem{gianazza2010}
D.~Gianazza, ``Forecasting workload and airspace configuration with neural
  networks and tree search methods,'' \emph{Artificial intelligence}, vol. 174,
  no. 7-8, pp. 530--549, 2010.

\bibitem{lee2008examining}
P.~Lee, J.~Mercer, B.~Gore, N.~Smith, K.~Lee, and R.~Hoffman, ``Examining
  airspace structural components and configuration practices for dynamic
  airspace configuration,'' in \emph{AIAA Guidance, Navigation and Control
  Conference and Exhibit}, 2008, p. 7228.

\bibitem{zelinski2011comparing}
S.~Zelinski and C.~F. Lai, ``Comparing methods for dynamic airspace
  configuration,'' in \emph{2011 IEEE/AIAA 30th Digital Avionics Systems
  Conference}.\hskip 1em plus 0.5em minus 0.4em\relax IEEE, 2011, pp. 3A1--1.

\bibitem{wang2021counter}
J.~Wang, Y.~Liu, and H.~Song, ``Counter-unmanned aircraft system (s)(c-uas):
  State of the art, challenges, and future trends,'' \emph{IEEE Aerospace and
  Electronic Systems Magazine}, vol.~36, no.~3, pp. 4--29, 2021.

\bibitem{yang2022multiagent}
Y.~Yang, J.~Yu, D.~Liu, S.-A. Lee, S.~Namilae, S.~Islam, H.~Gou, H.~Park, and
  H.~Song, ``Multiagent collaboration for emergency evacuation using
  reinforcement learning for transportation systems,'' \emph{IEEE Journal on
  Miniaturization for Air and Space Systems}, vol.~3, no.~4, pp. 232--241,
  2022.

\bibitem{yang2021machine}
Y.~Yang, K.~Zhang, H.~Song, and D.~Liu, ``Machine learning-enabled adaptive air
  traffic recommendation system for disaster evacuation,'' in \emph{2021
  IEEE/AIAA 40th Digital Avionics Systems Conference (DASC)}.\hskip 1em plus
  0.5em minus 0.4em\relax IEEE, 2021, pp. 1--8.

\bibitem{kopardekar2007initial}
P.~Kopardekar, K.~Bilimoria, and B.~Sridhar, ``Initial concepts for dynamic
  airspace configuration,'' in \emph{7th AIAA ATIO Conf, 2nd CEIAT Int'l Conf
  on Innov and Integr in Aero Sciences, 17th LTA Systems Tech Conf; followed by
  2nd TEOS Forum}, 2007, p. 7763.

\bibitem{kopardekar2008airspace}
P.~Kopardekar, K.~D. Bilimoria, and B.~Sridhar, ``Airspace configuration
  concepts for the next generation air transportation system,'' \emph{Air
  Traffic Control Quarterly}, vol.~16, no.~4, pp. 313--336, 2008.

\bibitem{Sergeevagenetic}
M.~Sergeeva, D.~Delahaye, C.~Mancel, and A.~Vidosavljevic, ``Dynamic airspace
  configuration by genetic algorithm,'' \emph{Journal of traffic and
  transportation engineering (English edition)}, vol.~4, no.~3, pp. 300--314,
  2017.

\bibitem{brinton2008airspace}
C.~R. Brinton and S.~Pledgie, ``Airspace partitioning using flight clustering
  and computational geometry,'' in \emph{2008 IEEE/AIAA 27th Digital Avionics
  Systems Conference}.\hskip 1em plus 0.5em minus 0.4em\relax IEEE, 2008, pp.
  3--B.

\bibitem{brinton2009airspace}
C.~Brinton, J.~Hinkey, and K.~Leiden, ``Airspace sectorization by dynamic
  density,'' in \emph{9th AIAA Aviation Technology, Integration, and Operations
  Conference (ATIO) and Aircraft Noise and Emissions Reduction Symposium
  (ANERS)}, 2009, p. 7102.

\bibitem{sabhnani2010flow}
G.~Sabhnani, A.~Yousefi, and J.~S. Mitchell, ``Flow conforming operational
  airspace sector design,'' in \emph{10th AIAA Aviation Technology,
  Integration, and Operations (ATIO) Conference}, 2010, p. 9377.

\bibitem{von2007tutorial}
U.~Von~Luxburg, ``A tutorial on spectral clustering,'' \emph{Statistics and
  computing}, vol.~17, pp. 395--416, 2007.

\bibitem{galluccio2012graph}
L.~Galluccio, O.~Michel, P.~Comon, and A.~O. Hero~III, ``Graph based k-means
  clustering,'' \emph{Signal Processing}, vol.~92, no.~9, pp. 1970--1984, 2012.

\bibitem{wang2016auto}
Y.~Wang, H.~Yao, and S.~Zhao, ``Auto-encoder based dimensionality reduction,''
  \emph{Neurocomputing}, vol. 184, pp. 232--242, 2016.

\bibitem{li2009spectral}
J.~Li, T.~Wang, I.~Hwang, and I.~Hwang, ``A spectral clustering based algorithm
  for dynamic airspace configuration,'' in \emph{9th AIAA Aviation Technology,
  Integration, and Operations Conference (ATIO) and Aircraft Noise and
  Emissions Reduction Symposium (ANERS)}, 2009, p. 7056.

\bibitem{gondara2016medical}
L.~Gondara, ``Medical image denoising using convolutional denoising
  autoencoders,'' in \emph{2016 IEEE 16th international conference on data
  mining workshops (ICDMW)}.\hskip 1em plus 0.5em minus 0.4em\relax IEEE, 2016,
  pp. 241--246.

\bibitem{Transportation_2017}
\BIBentryALTinterwordspacing
D.~o. Transportation, ``2015 flight delays and cancellations,'' Feb 2017.
  [Online]. Available:
  \url{https://www.kaggle.com/datasets/usdot/flight-delays}
\BIBentrySTDinterwordspacing

\bibitem{jiang2021spatial}
Y.~Jiang, S.~Niu, K.~Zhang, B.~Chen, C.~Xu, D.~Liu, and H.~Song,
  ``Spatial--temporal graph data mining for iot-enabled air mobility
  prediction,'' \emph{IEEE Internet of Things Journal}, vol.~9, no.~12, pp.
  9232--9240, 2021.

\bibitem{ng2001spectral}
A.~Ng, M.~Jordan, and Y.~Weiss, ``On spectral clustering: Analysis and an
  algorithm,'' \emph{Advances in neural information processing systems},
  vol.~14, 2001.

\end{thebibliography}
% Tell the compiler where to look for your references, IPCCC2022.bib here is the database that contains reference entries.

\end{document}